\documentclass[12pt]{article}
\usepackage{amsmath,amsfonts,latexsym,amsthm,amssymb}
\newcommand{\D}{D^\alpha}

\numberwithin{equation}{section}

\DeclareMathOperator{\const}{const}

\begin{document}

\newtheorem{lem}{Lemma}
\newtheorem{teo}{Theorem}
\newtheorem{cor}{Corollary}
\newtheorem{prop}{Proposition}

\pagestyle{plain}
\title{$L^p$ properties of non-Archimedean fractional differentiation operators}
\author{Anatoly N. Kochubei
\\ \footnotesize Institute of Mathematics,\\
\footnotesize National Academy of Sciences of Ukraine,\\
\footnotesize Tereshchenkivska 3, Kyiv, 01024 Ukraine\\
\footnotesize E-mail: \ kochubei@imath.kiev.ua}
\date{}
\maketitle

\vspace*{3cm}
\begin{abstract}
Let $\D$ be the Vladimirov-Taibleson fractional differentiation operator acting on \linebreak complex-valued functions on a non-Archimedean local field. The identity $\D D^{-\alpha}f=f$ was known only for the case where $f$ has a compact support. Following a result by Samko about the fractional Laplacian of real analysis, we extend the above identity in terms of $L^p$-convergence of truncated integrals. Differences between real and non-Archimedean cases are discussed.
\end{abstract}

\vspace{2cm}
{\bf Key words: }\ fractional differentiation operator; non-Archimedean local field; principal value operator; Riesz potentials

\medskip
{\bf MSC 2020}. Primary: 47G30. Secondary: 11S80, 35S05, 34A08.

{\bf Acknowledgement:} This work was funded in part under the research project "Markov evolutions in real and p-adic spaces" of the Dragomanov National Pedagogic University of Ukraine.

\newpage
\section{Introduction}

Contemporary mathematical analysis contains two parallel theories of nonlocal operators -- the fractional Laplacian of real analysis, its various interpretations and extensions (see the recent surveys \cite{KKK,Kw,St} as representatives of the vast literature on this subject) and the Vladimirov, or Vladimirov-Taibleson, fractional differentiation operator $\D$ acting on complex-valued functions on a non-Archimedean local field \cite{AKS,KKZ,K2001,VVZ}. Both theories originate from the notion of the Riesz potential (for the non-Archimedean case see \cite{Ta}), and the initial appearance of the fractional Laplacian was in the shape of a left inverse to the Riesz potential; see \cite{Sa2002} for further references.

Both the investigation of the fractional Laplacian and Vladimirov operator are far from their completion. For the latter, it is worth mentioning discoveries of such hidden structures as an isotropic Laplacian \cite{BGPW}, oscillation properties of the corresponding heat kernel \cite{BCW}, integral equations resembling the classical Volterra ones connected to the restriction of $\D$ to radial functions \cite{K2014,K2020_1,K2020_2}.

If $\D$, $\alpha >0$, is the Vladimirov operator, and $D^{-\alpha}$ is the Riesz potential (see Section 3 for the details), then $\D D^{-\alpha}\varphi =\varphi$ for a test function $\varphi$; such a function $\varphi$, in particular, must have a compact support. This important identity expressing the reciprocal inversion property of the operators $\D$ and $D^{-\alpha}$ can be extended to distributions with compact supports, but the last assumption cannot be dropped (\cite{VVZ}, Chapter 2.IX).

In this paper, we prove a kind of the above inversion property for $\varphi \in L^p$, with the hypersingular integral operator $\D$ understood in the sense of principal value, with different meaning of convergence depending on the vaule of $\alpha$. Similar results for the fractional Laplacian were obtained by Samko \cite{Sa76,Sa2002}. However there is an interesting difference between the real and non-Archimedean cases.

The truncated operator $D^\alpha_\varepsilon (D^{-\alpha }\varphi)$ (with a removed singularity) is expressed via $\varphi$  through integration with an averaging kernel $\mathcal K$. In the real case, $\mathcal K$ has a power-like decay at infinity. The non-Archimedean phenomenon is that $\mathcal K$ has a compact support.

Another non-Archimedean feature is the sufficiency to perform analytic calculations only in the one-dimensional case. It is shown that the transition to the multi-dimensional situation is equivalent to considering the whole setting over a wider field, an unramified extension of the initial one.

The structure of this paper is as follows. In Section 2, we present basic information about the Vladimirov operator and field extensions. In Section 3, we expose our main results for one-dimensional operators, while in Section 4 we explain the interpretation of the multi-dimensional operator as a one-dimensional one over an unramified extension.

\section{Preliminaries}

{\bf 2.1. Local fields.} A non-Archimedean local field is a non-discrete totally disconnected locally compact topological field. Such a field $K$ is isomorphic either to a finite extension of the field $\mathbb Q_l$ of $l$-adic numbers (where $l$ is a prime number), if $K$ has characteristic 0, or to the field of formal Laurent series with coefficients from a finite field, if $K$ has
a positive characteristic. A summary of main notions and results regarding local fields is given, for example, in \cite{K2001}.

Any local field is endowed with an absolute value $|\cdot |_K$ possessing the following properties:
1) $|x|_K=0$ if and only if $x=0$, 2) $|xy|_K=|x|_K\cdot |y|_K$,
3) $|x+y|_K\le \max (|x|_K,|y|_K)$. The last property called the ultrametric one implies that $|x+y|_K=|x|_K$, if $|y|_K<|x|_K$.

Denote $O=\{ x\in K:\ |x|_K\le 1\}$, $P=\{ x\in K:\ |x|_K<1\}$. $O$ is a subring of $K$, and $P$ is an ideal in $O$ containing such an element $\beta$ that $P=\beta O$. The quotient ring $\bar{K}=O/P$ is
a finite field called the residue field. The absolute value is called normalized, if $|\beta |_K=q^{-1}$ where $q$ is the cardinality of $O/P$. Unless stated otherwise, the absolute value used below is assumed normalized.
The normalized absolute value takes the values $q^N$, $N\in \mathbb Z$. For $K=\mathbb Q_l$ we have $\beta =l$ and $q=l$ where $l$ is understood in the first equality as a member of $\mathbb Q_l$, and in the second equality as a natural number.

The additive group of a local field is self-dual, so that the Fourier analysis on $K$ resembles the classical one. Let $\chi$ be a fixed non-constant complex-valued additive character on
$K$. Then any other additive character has the form $x\mapsto \chi (ax)$ for some $a\in K$. Below we assume that $\chi$ is a rank
zero character, that is $\chi (x)\equiv 1$ for $x\in O$, while $\chi
(x_0)\ne 1$ for some $x_0\in K$ with $|x_0|_K=q$.

Denoting by $dx$ the Haar measure on the
additive group of $K$ (normalized in such a way that the measure
of $O$ equals 1) we define the Fourier transform of a complex-valued function $f\in L^1(K)$ as
$$
(\mathcal F)(\xi )=\widetilde{f}(\xi )=\int\limits_K\chi (x\xi )f(x)\,dx,\quad \xi
\in K.
$$
If $\mathcal F f=\widetilde{f}\in L^1(K)$, then we
have the inversion rule
$$
f(x)=\int\limits_K\chi (-x\xi )\widetilde{f}(\xi )\,d\xi .
$$

Below we use standard integration formulas \cite{K2001,Vl,VVZ}. Here we list some of them; for $n\in \mathbb Z$, $\alpha >0$,
\begin{equation*}
\int\limits_{|x|_K\le q^n}dx=q^n;\quad \int\limits_{|x|_K=q^n}dx=\left( 1-\frac1q \right)q^n.
\end{equation*}

\begin{equation}
\label{2.1}
\int\limits_{|x|_K\le q^n}|x|_K^{\alpha -1}\,dx=\frac{1-q^{-1}}{1-q^{-\alpha }}q^{\alpha n}.
\end{equation}

\begin{equation}
\label{2.2}
\int\limits_{|x|_K=q^n}|x-a|_K^{\alpha -1}\,dx=\frac{q-2+q^{-\alpha}}{q(1-q^{-\alpha})}|a|_K^\alpha,\ \text{if $|a|_K=q^n$};
\end{equation}

\begin{equation}
\label{2.3}
\int\limits_{|x|_K\le q^n}\log |x|_K\,dx=\left( n-\frac1{q-1}\right) q^n\log q;
\end{equation}

\begin{equation}
\label{2.4}
\int\limits_{|x|_K=q^n}\log |x-a|_K\,dx=\left[ \left( 1-\frac1q\right)\log |a|_K-\frac{\log q}{q-1}\right] |a|_K,\ \text{if $|a|_K=q^n$}.
\end{equation}

A function $f:\ K\to \mathbb C$ is called locally constant, if there exists such an integer $k$ that for any $x\in K$
$$
f(x+x')=f(x), \quad \text{whenever $|x'|\le q^{-k}$}.
$$

The vector space $\mathcal D(K)$ of all locally constant functions with compact supports is called the Bruhat-Schwartz space of test functions. The Fourier transform preserves $\mathcal D(K)$. In some cases, it is convenient to use a more narrow Lizorkin space
$$
\Phi (K)=\left\{ \varphi \in \mathcal D(K):\ \int\limits_K \varphi (x)\,dx=0\right\};
$$
for the details see \cite{AKS}.

The above Fourier analysis is extended easily to functions on $K^n$. The natural non-Archimedean norm on $K^n$ is
\begin{equation}
\label{2.5}
\| (x_1,\ldots ,x_n)\|=\max\limits_{1\le j\le n}|x_j|_K.
\end{equation}

\medskip
{\bf 2.2. Field extensions.}  If a local field $K$ is a subfield of a local field $L$, then $L$ is called an extension of $K$ (which is denoted $L/K$). The extension $L$ can be considered as a vector space over $K$. $L$ is called a finite extension, if it is a finite-dimensional vector space over $K$. Its dimension is called the degree of the extension. Any basis of $L$ over $K$ is called a basis of the extension.

An operator of multiplication in $L$ by an element $\xi$ can be considered as a linear operator in the $K$-vector space. If the linear function $\xi \mapsto \operatorname{Tr}(\xi)$ does not vanish identically, then the extension is called separable. Such are all finite extensions of a field of characteristic 0. The above notion of separability makes sense also for finite fields $\bar{K},\bar{L}$.

There exists a detailed theory of extensions of local fields; see \cite{FV,Weil,Weiss}. In this paper, we need only one class, the unramified extensions.

A finite extension $L$ of a local field $K$ is called unramified, if $\bar{L}/\bar{K}$ is a separable extension of the same degree as $L/K$. Any local field $K$ has a unique (up to isomorphism) unramified extension of each given degree $n\ge 1$. If $\beta$ is a prime element of $K$, it remains a prime element of each unramified extension. If $L$ is the unramified extension of $K$ of a degree $n$, then the cardinality of the residue field equals $q^n$ where $q$ is the cardinality of the residue field . Each unramified extension is separable.

The unramified extension $L$ of degree $n$ has, as a vector space over $K$, a canonical basis consisting of representatives of a basis in $\bar{L}$ over $\bar{K}$. If $x\in L$ has the coefficients $x_1,\ldots ,x_n\in K$ of the expansion with respect to the canonical basis, then the normalized absolute value $|x|_L$ has the representation \cite{Ta76}
\begin{equation}
\label{2.6}
|x|_L=\left( \max\limits_{1\le j\le n}|x_j|_K\right)^n.
\end{equation}

Along with the normalized absolute value $|\cdot |_L$, there is another absolute value on $L$ extending the absolute value from $K$, namely $\|x\|_L=|x|_L^{1/n}$. Comparing (\ref{2.6}) and (\ref{2.5}) we see that the above canonical basis is an orthonormal basis (see \cite{Sch}) in $L$ as a Banach space with the norm $\|\cdot \|_L$. Therefore the canonical basis defines an isometric isomorphism between $L$ and $K^n$.

In particular, this isomorphism defines an invariant measure on the additive group $K^n$, and the measure of the unit ball in $K^n$ equals 1. By the uniqueness of the Haar measure, we have $dx=dx_1\cdots dx_n$.

\medskip
\section{The Vladimirov operator}.

{\bf 3.1. Definitions.} On a test function $\varphi \in \mathcal D(K)$, the fractional differentiation operator $\D$, $\alpha >0$, is defined as follows:
$$
(\D \varphi )(x)=\mathcal F^{-1}\left[ |\xi |_K^\alpha (\mathcal F(\varphi))(\xi )\right] (x), \quad x\in K.
$$
While $D^\alpha$ does not preserve $\mathcal D(K)$, it preserves the Lizorkin space $\Phi (K)$ \cite{AKS}.

The operator $\D$ admits a hypersingular integral representation \cite{K2001,VVZ}
\begin{equation}
\label{3.1}
\left( D^\alpha \varphi \right) (x)=c_\alpha \int\limits_K |y|_K^{-\alpha -1}[\varphi (x-y)-\varphi (x)]\,dy,\quad c_\alpha =\frac{1-q^\alpha}{1-q^{-\alpha -1}}.
\end{equation}
The expression (\ref{3.1}) makes sense for wider classes of functions.

A right inverse to $\D$ defined on $\mathcal D(K)$ is given by the Riesz potential $D^{-\alpha}$ \cite{VVZ}
$$
(D^{-\alpha }\varphi )(x)=d_\alpha \int\limits_K|x-y|_K^{\alpha -1}\varphi (y)\,dy,\quad \alpha \ne 1,
$$
where $d_\alpha =\dfrac{1-q^{-\alpha}}{1-q^{\alpha -1}}$. For $\alpha =1$, we have a right inverse $D^{-1}$ on $\Phi (K)$:
$$
(D^{-1}\varphi )(x)=d_1 \int\limits_K\log |x-y|_K\varphi (y)\,dy,\quad d_1=\frac{1-q}{q\log q}.
$$

As we emphasized in Introduction, our task is to extend the above definitions in such a way that the identity $\D D^{-\alpha}\varphi =\varphi$ remains valid for ``worse'' functions $\varphi$. We will understand $\D$ as a principal value integral
\begin{equation}
\label{3.2}
\left( D^\alpha u \right) (x)=\lim\limits_{\varepsilon\to 0}\left( D^\alpha_\varepsilon u \right) (x)
\end{equation}
where
\begin{equation}
\label{3.3}
\left( D^\alpha_\varepsilon u \right) (x)=c_\alpha \int \limits_{|y-x|_K\ge \varepsilon}\frac{u(y)-u(x)}{|y-x|_K^{\alpha +1}}\,dy,
\end{equation}
and the meaning of the limit will be specified later.

\medskip
{\bf 3.2. The case where $0<\alpha <1$.} Suppose that
\begin{equation}
\label{3.4}
u= D^{-\alpha}\varphi \ (0<\alpha <1),\quad \varphi \in L^p(K),1\le p<\frac1\alpha.
\end{equation}

It is known (see \cite{Ta}, Chapter III, Theorem 4.9) that under the assumption (\ref{3.4}) $D^{-\alpha}\varphi$ exists a.e., and if $p>1$, then $D^{-\alpha}\varphi\in L^r(K)$ where $0<\frac1r=\frac1p -\alpha$. If $p=1$, then $D^{-\alpha}\varphi$ belongs to the weak space $L^{r,\infty}(K)$ (for the definition and properties in the analysis on general groups see \cite{Gr}). The existence of the integral in (\ref{3.3}) follows from the Young inequalities for strong and weak $L^p$-spaces; see Theorems 1.2.12 and 1.2.13 in \cite{Gr}.

It is convenient to set $\varepsilon =q^{-\nu}$, $\nu \in \mathbb N$. Choose $\sigma \in K$ in such a way that $|\sigma|_K=\varepsilon$.

\medskip
\begin{prop}
Under the assumptions (\ref{3.4}), the truncated operator $D^\alpha_\varepsilon u$, $u=D^{-\alpha }\varphi$, admits a representation
\begin{equation}
\label{3.5}
\left( D^\alpha_\varepsilon u \right) (x)=c_\alpha d_\alpha \int\limits_KR(\tau)\varphi (x-\sigma \tau )\,d\tau
\end{equation}
where
\begin{equation}
\label{3.6}
R(\tau )=\begin{cases}
1-\left( 1-\frac1q\right) \frac1{1-q^{-\alpha}}|\tau |_K^{\alpha-1}, & \text{if $|\tau|_K<1$},\\
0, & \text{if $|\tau|_K\ge 1$}.\end{cases}
\end{equation}
For the kernel $R_1(\tau )=c_\alpha d_\alpha R(\tau)$, we have
\begin{equation}
\label{3.7}
R_1(\tau )>0, \text{as $|\tau|_K<1$};\quad \int\limits_K R_1(\tau )\,d\tau =1.
\end{equation}
\end{prop}

\medskip
{\it Proof.} Substituting $u=D^{-\alpha }\varphi$ into (\ref{3.3}) we have an absolutely convergent double integral, in which we can change the order of integration getting
$$
\left( D^\alpha_\varepsilon u \right) (x)=c_\alpha d_\alpha \int\limits_K \varphi (x-t)\,dt \int\limits_{|s|_K\ge \varepsilon}|s|_K^{-\alpha -1}\left( |s+t|_K^{\alpha -1}-|t|_K^{\alpha -1}\right)\,ds.
$$
After the change of variables $t=\sigma \tau$, $s=\sigma \xi$, we come to the expression (\ref{3.5}) with
\begin{equation}
\label{3.8}
R(\tau )=\int\limits_{|\xi |_K\ge 1}|\xi |_K^{-\alpha -1}\left( |\xi +\tau|_K^{\alpha -1}-|\tau|^{\alpha -1}\right) \,d\xi.
\end{equation}

Let us calculate $I_1=\int\limits_{|\xi |_K\ge 1}|\xi |_K^{-\alpha -1}|\xi +\tau|_K^{\alpha -1}\,d\xi $. We have to consider three cases.

1) $|\tau|_K<1$. Then $|\xi +\tau|_K=|\xi|_K$, so that
$$
I_1=\int\limits_{|\xi |_K\ge 1}|\xi |_K^{-2}\,d\xi =\sum\limits_{j=0}^\infty \int\limits_{|\xi |_K=q^j}|\xi |_K^{-2}\,d\xi=\left( 1-\frac1q\right)\sum\limits_{j=0}^\infty q^{-j}=1.
$$

2) $|\tau|_K=1$. After the change $\xi =\tau r$ we obtain
\begin{equation}
\label{3.9}
I_1=|\tau|_K^{-1}\int \limits_{|r|_K\ge |\tau|_K^{-1}}|r|_K^{-\alpha -1}|1+r|_K^{\alpha -1}\,dr.
\end{equation}

Below we use (\ref{3.9}) also for $|\tau|_K>1$. Meanwhile, for $|\tau|_K=1$ it follows from (\ref{3.9}) that
$$
I_1=\int\limits_{|r|_K\ge 1}|r|_K^{-\alpha -1}|1+r|_K^{\alpha -1}\,dr=\int\limits_{|r|_K=1}|1+r|_K^{\alpha -1}\,dr+\int\limits_{|r|_K>1}|r|_K^{-2}\,dr,
$$
and by (\ref{2.2}),
$$
I_1=\frac{q-2+q^{-\alpha}}{q(1-q^{-\alpha})}+\left( \frac1q\right)\sum\limits_{j=1}^\infty q^{-j}=\frac{q-1}{q(1-q^{-\alpha})}.
$$

3) $|\tau|_K>1$. Let $|\tau|_K=q^N$, $N\ge 1$. Proceeding from (\ref{3.9}) and using (\ref{2.2}) we find that
\begin{multline*}
I_1=|\tau|_K^{-1}\sum\limits_{j=-N}^\infty q^{-j(\alpha +1)}\int\limits_{|r|_K=q^j}|1+r|_K^{\alpha -1}\,dr\\
=|\tau|_K^{-1}\left\{ \left( \frac1q\right)\sum\limits_{j=-N}^{-1}q^{-\alpha j}+\int\limits_{|r|_K=1}|1+r|_K^{\alpha -1}\,dr+\left( \frac1q\right)\sum\limits_{j=1}^\infty q^{-j}\right\} \\
=|\tau|_K^{-1}\left\{ \left( \frac1q\right)\frac{q^{\alpha (N+1)}-q^\alpha}{q^\alpha -1}+ \frac{q-2+q^{-\alpha}}{q(1-q^{-\alpha})}+\frac1q\right\} \\
=\left(1-\frac1q\right)\frac1{1-q^{-\alpha}}|\tau|_K^{\alpha -1}+|\tau|_K^{-1}\left[ -\frac{q-1}{q(1-q^{-\alpha})}+ \frac{q-2+q^{-\alpha}}{q(1-q^{-\alpha})}+\frac1q\right] \\ =\left(1-\frac1q\right)\frac1{1-q^{-\alpha}}|\tau|_K^{\alpha -1},
\end{multline*}
since the expression in brackets equals zero.

In order to compute $R(\tau )$, we need to calculate also
$$
I_2=|\tau|_K^{\alpha -1}\int\limits_{|\xi|_K\ge 1}|\xi|_K^{-\alpha -1}\,d\xi.
$$
We have
$$
I_2=\left(1-\frac1q\right)|\tau|_K^{\alpha -1}\sum\limits_{j=0}^\infty q^{-\alpha j}=\left(1-\frac1q\right)\frac1{1-q^{-\alpha}}|\tau|_K^{\alpha -1}.
$$

Comparing this with the expressions for $I_1$ for the above three cases, we come to (\ref{3.6}).

Using (\ref{2.1}) and (\ref{3.6}), we find that
$$
\int\limits_{|\tau|_K\le q^{-1}}R(\tau )\,d\tau =q^{-1}-\frac{(1-q^{-1})^2q^{-\alpha}}{(1-q^{-\alpha})^2}=\frac{q^{-1}+q^{-2\alpha -1}-q^{-\alpha}-q^{-2-\alpha}}{(1-q^{-\alpha})^2}.
$$
On the other hand,
$$
c_\alpha d_\alpha=\frac{(1-q^{-\alpha})^2}{q^{-2\alpha -1}-q^{-\alpha}-q^{-\alpha -2}+q^{-1}},
$$
so that $\int\limits_K R_1(\tau )\,d\tau =1$.

If $|\tau|_K=q^{-j}$, $j\ge 1$, then $R(\tau)$ is a decreasing function of $j$. Therefore if we show that $R(\tau )<0$ for $j=1$, then we will prove that $R(\tau )<0$ for $|\tau|_K<1$. In fact, for $|\tau|_K=q^{-1}$ we have
$$
1-\frac{q-1}q\frac1{1-q^{-\alpha}}q^{1-\alpha}=\frac{1-q^{1-\alpha}}{1-q^{-\alpha}}<0.
$$
Since $c_\alpha <0$, $d_\alpha >0$, we find that $R_1(\tau )>0$. $\qquad \blacksquare$

Now we can prove that the principal value operator on $L^p(K)$ is indeed the left inverse of the Riesz potential.

\medskip
\begin{teo}
Under the assumptions (3.4),
$$
\| D_\varepsilon^\alpha (D^{-\alpha} \varphi) -\varphi \|_{L^p(K)}\longrightarrow 0, \text{ as $\varepsilon =q^{-\nu}\to 0$}.
$$
\end{teo}

\medskip
{\it Proof}. We have proved that
$$
\left( D_\varepsilon^\alpha (D^{-\alpha} \varphi )\right)(x)-\varphi (x)=\int\limits_K R_1(\tau)[\varphi (x-\sigma \tau)-\varphi (x)]\,dx.
$$
By Minkowski's inequality, we get
\begin{equation}
\label{3.10}
\| D_\varepsilon^\alpha (D^{-\alpha} \varphi) -\varphi \|_{L^p(K)}\le \int\limits_K R_1(\tau)\omega_p(\varphi ,\sigma \tau)\,d\tau
\end{equation}
where $\omega_p (\varphi, h)=\|\varphi (\cdot)-\varphi (\cdot -h)\|_{L^p(K)}$. By the dominated convergence theorem, the left-hand side of (\ref{3.10}) tends to zero. $\qquad \blacksquare$

{\bf Remark.} Proposition 1 remains valid for $\alpha >1$, if $\varphi \in \Phi (K)$. We have the representation (\ref{3.5}) with the kernel (\ref{3.6}). Of the properties (\ref{3.7}), we cannot assert the positivity of the kernel $R_1$; the normalization identity remains valid.

\medskip
{\bf 3.3. The case where $\alpha =1$.} First we prove for $\alpha =1$ an analog of Proposition 1. Our assumptions for this case are of a different nature:
\begin{equation}
\label{3.11}
\varphi \in L^\infty_{\text loc}(K),\quad \varphi =O(|t|_K^{-\beta }), \beta >1,t\to \infty.
\end{equation}

\medskip
\begin{prop}
Under the assumptions (\ref{3.11}), the truncated operator $D^1_\varepsilon u$, $u=D^{-1}\varphi$, admits a representation (\ref{3.5}) with $\alpha =1$,
$$
R(\tau )=\begin{cases}
\frac{\log q}{q-1}-\log |\tau|_K, & \text{if $|\tau|_K<1$},\\
0, & \text{if $|\tau|_K\ge 1$},\end{cases}
$$
where the kernel $R_1(\tau)=c_1d_1R(\tau)$ possesses the properties (\ref{3.7}).
\end{prop}

\medskip
{\it Proof}. In the present case
$$
\left( D_\varepsilon^1 u\right) (x)=c_1\int\limits_{|z-x|_K\ge \varepsilon}\frac{u(z)-u(x)}{|z-x|^2}\,dz,\quad c_1=-\frac{q^2}{q+1},
$$
$$
\left( D^{-1}\varphi \right) (x)=d_1\int\limits_K \log |y-x|_K\varphi (y)\,dy,\quad d_1=\frac{q-1}{q\log q},
$$
so that for $u= D^{-1}\varphi$ we have
$$
\left( D_\varepsilon^1 u\right) (x)=c_1d_1\int\limits_{|z-x|_K\ge \varepsilon}|z-x|_K^{-2}\,dz\int\limits_K(\log |y-z|_K-\log |y-x|_K)\varphi (y)\,dy.
$$

After the changes of variables $x-y=t, z-x=s$ we get
\begin{equation}
\label{3.12}
\left( D_\varepsilon^1 u\right) (x)=c_1d_1\int\limits_{|s|_K\ge \varepsilon}|s|_K^{-2}\,ds\int\limits_K(\log |s+t|_K-\log |t|_K)\varphi (x-t)\,dt.
\end{equation}
The integral in (\ref{3.12}) is absolutely convergent. To prove that, consider, for a fixed $x\in K$, the integral
$$
\lambda_x(s)=\int\limits_K (\log |s+t|_K-\log |t|_K)\varphi (x-t)\,dt=\int\limits_{|t|_K\le |s|_K} (\log |s+t|_K-\log |t|_K)\varphi (x-t)\,dt
$$
(the last equality follows from the ultrametric property).

Let $|s|_K>|x|_K$. If $|t|_K>|x|_K$, then $|x-t|_K=|t|_K$, so that
$$
|\varphi (x-t)|\le C|t|_K^{-\beta},\quad \text{as $|t|_K>|x|_K$}.
$$

Let us write $\lambda_x(s)=\lambda_x^{(1)}(s)+\lambda_x^{(2)}(s)$ where the summands correspond to the integration on $\{ t:\ |t|_K\le |x|_K\}$ and $\{ t:\ |x|_K<|t|_K\le |s|_K\}$ respectively.

If $|t|_K\le |x|_K$, then $|x-t|_K\le |x|_K$, so that $|\varphi (x-t)|\le \const$ and
$$
\left| \lambda_x^{(1)}(s)\right| \le C\int\limits_{|t|_K\le |x|_K} (\log |s|_K+\log |t|_K)\, dt\le C|\log |s|_K|.
$$

Next, using (2.1) we find that
\begin{multline*}
\left| \lambda_x^{(2)}(s)\right| \le C\int\limits_{|x|_K<|t|_K\le |s|_K} (\log |s+t|_K+\log |t|_K)|t|_K^{-\beta}\, dt\\
\le C_1|\log |s|_K|\int\limits_{|t|_K\le |s|_K}|t|_K^{-\beta}\, dt=C_2|\log |s|_K|\cdot |s|_K^{1-\beta}.
\end{multline*}

Therefore the integral (\ref{3.12}) is absolutely convergent, and we may change the order of integration
$$
\left( D_\varepsilon^1 u\right) (x)=c_1d_1\int\limits_K R(\tau)\varphi (x-\sigma \tau)\,d\tau
$$
where $|\sigma|_K=\varepsilon$, and after easy transformations,
$$
R(\tau)=\int\limits_{|\xi |_K\ge 1}|\xi |_K^{-2}(\log |\xi +\tau|_K-\log |\tau|_K)\,d\xi=J_1-J_2
$$
where
\begin{equation}
\label{3.13}
J_1=\int\limits_{|\xi |_K\ge 1}|\xi |_K^{-2}\log |\xi +\tau|_K\,d\xi,
\end{equation}
\begin{equation}
\label{3.14}
J_2=\log |\tau|_K \int\limits_{|\xi |_K\ge 1}|\xi|_K^{-2}\,d\xi.
\end{equation}

The change of variables $\xi =\tau r$ brings $J_1$ to the form
\begin{equation}
\label{3.15}
J_1=|\tau|_K^{-1}\int\limits_{|r|_K\ge |\tau|_K^{-1}}|r|_K^{-2}(\log |\tau|_K+\log |1+r|_K\,dr.
\end{equation}

If $|\tau|_K<1$, then by (\ref{3.13}) and the formula 0.231.2 from \cite{GR} we get
$$
J_1=\int\limits_{|\xi |_K\ge 1}|\xi |_K^{-2}\log |\xi|_K\,d\xi=\left(1-\frac1q\right)\log q\sum\limits_{j=0}^\infty jq^{-j}=\frac{\log q}{q-1},
$$
so that
\begin{equation}
\label{3.16}
J_1=\frac{\log q}{q-1},\quad |\tau|_K<1.
\end{equation}

If $|\tau|_K=1$, then we use (\ref{3.15}), (\ref{2.4}) and the same formula from \cite{GR}:
\begin{multline*}
J_1=\int\limits_{|r|_K\ge 1}|r|_K^{-2}\log |1+r|_K\,dr=\int\limits_{|r|_K=1}\log |1+r|_K\,dr\\
+\sum\limits_{j=1}^\infty \int\limits_{|r|_K=q^j}|r|_K^{-2}\log |r|_K\,dr=-\frac{\log q}{q-1}+\left(1-\frac1q\right)\log q\sum\limits_{j=1}^\infty jq^{-j}=0
\end{multline*}
so that
\begin{equation}
\label{3.17}
J_1=0,\quad |\tau|_K=1.
\end{equation}

Consider the case where $|\tau|_K>1$. Let us write $|\tau|_K=q^N$, $N\ge 1$. Using the representation (\ref{3.15}), the integration formula (\ref{2.4}), the formula 0.231.2 from \cite{GR}, and noticing that $\log |1+r|_K=0$ as $|r|_K<1$, we get
\begin{multline*}
J_1=q^{-N}\int\limits_{|r|_K\ge q^{-N}}|r|_K^{-2}(N\log q+\log |1+r|_K)\,dr\\
=Nq^{-N}\log q\int\limits_{|r|_K\ge q^{-N}}|r|_K^{-2}\,dr+q^{-N}\int\limits_{|r|_K\ge q^{-N}}|r|_K^{-2}\log |1+r|_K\,dr\\
=Nq^{-N}\left(1-\frac1q\right)\log q\sum\limits_{j=-N}^\infty q^{-j}+q^{-N}\int\limits_{|r|_K=1}\log |1+r|_K\,dr+q^{-N}\sum\limits_{j=1}^\infty q^{-2j}\int\limits_{|r|_K=q^j}\log |r|_Kdr\\
=N\log q-q^{-N}\frac{\log q}{q-1}+q^{-N}\left(1-\frac1q\right)\log q\sum\limits_{j=1}^\infty jq^{-j}=N\log q=\log |\tau|_K,
\end{multline*}
so that
\begin{equation}
\label{3.18}
J_1=\log |\tau|_K,\quad |\tau|_K>1.
\end{equation}

On the other hand, by (\ref{3.14}),
$$
J_2=\log |\tau|_K\cdot \left(1-\frac1q\right) \sum\limits_{j=0}^\infty q^{-j}=\log |\tau|_K.
$$
Comparing this with (\ref{3.16}), (\ref{3.17}) and (\ref{3.18}),  we obtain the required representation formulated in Proposition 2.

Note that $c_1d_1=-\dfrac{q(q-1)}{(q+1)\log q}$, thus for $|\tau|_K<1$,
$$
R_1(\tau)=\frac{q}{q+1}-\frac{q(q-1)}{(q+1)\log q}\log |\tau|_K,
$$
and we can calculate the integral using (\ref{2.3}):
$$
\int\limits_{|\tau|_K<1}R_1(\tau)\,d\tau=\frac1{(q+1)}-\frac{q(q-1)}{(q+1)\log q}\left( -1-\frac1{q-1}\right))q^{-1}\log q,
$$
and it is easy to find that $\int\limits_{|\tau|_K<1}R_1(\tau)\,d\tau=1$. It is also clear that $R_1(\tau)>0$, as $|\tau|_K<1$. $\qquad\blacksquare$

As in the previous case, we have the following result about the principal value operator $D^1$.

\medskip
\begin{teo}
Under the assumptions (\ref{3.11}),
$$
\left\|D _\varepsilon^1(D^{-1}\varphi )-\varphi\right\|_{L^1(K)}\longrightarrow 0,\quad \text{as $\varepsilon=q^{-\nu}\to 0$}.
$$
\end{teo}

\medskip
The {\it proof} based on Proposition 2 is similar to that of Theorem 1. $\qquad\blacksquare$

\medskip\section{Multi-Dimensional Operators}

Let us consider the multi-dimensional operator
\begin{equation}
\label{4.1}
\left( D^\alpha_{K^n}v\right) (x)=\frac{1-q^\alpha}{1-q^{-\alpha -n}}\int\limits_{K^n}\frac{v(z)-v(x)}{\|z-x\|^{n+\alpha}}\,dz,\quad x\in K^n
\end{equation}
(called the Taibleson operator in \cite{AKS}). Here
$$
\| x\|=\max\limits_{1\le j\le n}|x_j|_K,\quad x=(x_1,\ldots ,x_n).
$$

Let $L$ be an unramified extension of $K$ of degree $n$. As explained in Section 2.2, the expansion with respect to a canonical basis in $L$ defines an isometric linear isomorphism between $L$ and $K^n$, which identifies $D^\alpha_{K^n}$ defined in (\ref{4.1}) with the operator
\begin{equation}
\label{4.2}
\left( D_L^\gamma u\right) (x)=\frac{1-q^{n\gamma}}{1-q^{-n(\gamma +1)}}\int\limits_L\frac{u(z)-u(x)}{|z-x|_L^{\gamma +1}}\,dz
\end{equation}
where $\gamma =\dfrac{\alpha}n$.

Applying the one-dimensional results to the operator (\ref{4.2}) we prove the $L^p$- convergence of the truncated operator for $0<\alpha <n$, $1\le p<\frac{n}\alpha$. The kernel $R$ (compare (\ref{3.6}))  equals
$$
R(\tau )=\begin{cases}
1-\left(1-\frac1{q^n}\right)\frac1{1-q^{-\alpha n}}\|\tau \|^{\frac{\alpha}n-1}, & \text{if $\|\tau\|<1$},\\
0, & \text{if $\|\tau\|\ge 1$}.\end{cases}
$$

The case $\alpha =n$ is reduced similarly to Proposition 2 and Theorem 2. In the conditions (\ref{3.11}) we need $\beta >n$ now. In the expression for the kernel $R$, we have to substitute $q^n$ for $q$.

\end{document}